\documentclass{article}
\usepackage{amsmath}
\usepackage{amssymb}
\usepackage{mathrsfs}
\usepackage{amsfonts}
\usepackage{multirow}
\usepackage{cases}
\allowdisplaybreaks[4]

\newtheorem{Theorem}{Theorem}[part]
\newtheorem{Definition}{Definition}[part]
\newtheorem{Proposition}{Proposition}[part]

\newtheorem{Lemma}{Lemma}[part]

\newtheorem{Algorithm}{Algorithm}[part]
\newtheorem{Condition}{Condition}[part]

\def \ep{\hbox{ }\hfill$\Box$}

\begin{document}
\title{An algorithm to find the spectral radius of nonnegative tensors and its convergence analysis
\thanks{This work was supported by the National Natural Science Foundation of China(Grant No. 10871105) and Scientific Research Foundation for the Returned Overseas Chinese Scholars, State Education Ministry. }}

\author{Yuning Yang\\
\small{School of Mathematics Science and LPMC }\\
\small{Nankai University}\\
\small{Tianjin 300071, P.R. China}\\
\small{Email:nk0310145@gmail.com }\vspace{4mm} \\
Qingzhi Yang\\
\small{School of Mathematics Science and LPMC }\\
\small{Nankai University}\\
\small{Tianjin 300071, P.R. China}\\
\small{Email:qz-yang@nankai.edu.cn }\vspace{4mm} \\
Yanguang Li\\
\small{School of Mathematics Science and LPMC }\\
\small{Nankai University}\\
\small{Tianjin 300071, P.R. China}\\
\small{Email:lyg2004hot@163.com }\vspace{4mm} \\
}

\maketitle

\begin{abstract}
In this paper we propose an iterative algorithm to find out the
spectral radius of nonnegative tensors. This algorithm is an
extension of the smoothing method for finding the largest eigenvalue
of a nonnegative matrix \cite{s14}. For nonnegative irreducible
tensors, we establish the converges of the algorithm. Finally we
report some numerical results and conclude this paper with some
remarks.\vspace{3mm}

\noindent {\bf Keywords } nonnegative tensor, spectral radius, smoothing method, diagonal transformation\\
{\bf MSC }  74B99, 15A18, 15A69

\hspace{2mm}\vspace{2mm}

\end{abstract}

\section{Introduction}
\setcounter{equation}{0} \setcounter{Assumption}{0}
\setcounter{Theorem}{0} \setcounter{Proposition}{0}
\setcounter{Corollary}{0} \setcounter{Lemma}{0}
\setcounter{Definition}{0} \setcounter{Remark}{0}
\setcounter{Algorithm}{0}  %\hspace*{\parindent}

Eigenvalue problems of high order tensor have become an important
topic of study in a new applied mathematics branch, numerical
multilinear algebra, and they have a wide range of practical
applications, for more references, see \cite{s2,s5,s6,s9,s10,s13}.
In recent years, the largest eigenvalue problem for nonnegative
tensors has attracted special attention. Chang \emph{et al}
\cite{s1} generalized the Perron-Frobenius theorem from nonnegative
irreducible matrices to nonnegative irreducible tensors. Ng
\emph{et al} \cite{s12} gave a method to find the largest eigenvalue
of a nonnegative irreducible tensor. Yang and Yang \cite{s11}
defined the spectral radius of a tensor and gave further results for the 
Perron-Frobenius theorem and proved that the spectral radius is the
largest eigenvalue of any nonnegative tensor and all eigenvalues
with the spectral radius as their modulus distribute uniformly on
the circle. In this paper, we propose a method to find the spectral
radius of a class of nonnegative tensors. This method is an
extension of a method in \cite{s14} for calculating the spectral
radius of a nonnegative matrix . We show that for  a nonnegative
irreducible tensor, the sequence generated by this algorithm
converges to the spectral radius.

This paper is organized as follows: In section 2 we recall some
definitions and theorems; we give our algorithm in section 3 and lay
down the proof of the algorithm in section 4; some numerical results
are reported in section 5.

We first add a comment on the notation that is used in the sequel. Vectors are
written as lowercase letters $(x,y,\ldots)$, matrices correspond to
italic capitals $(A,B,\ldots)$, and tensors are written as
calligraphic capitals $(\mathcal{A}, \mathcal{B}, \cdots)$. The
entry with row index $i$ and column index $j$ in a matrix $A$, i.e.
$(A)_{ij}$ is symbolized by $a_{ij}$(also $(\mathcal{A})_{i_1\cdots
i_p,j_1\cdots j_q} = a_{i_1\cdots i_p,j_1\cdots j_q})$. The symbol
$|\cdot|$ used on a matrix $A$(or tensor $\mathcal{A}$) means that
$(|A|)_{ij} = |a_{ij}|$(or $(|\mathcal{A}|)_{i_1\cdots i_p,j_1\cdots
j_q} = |a_{i_1\cdots i_p,j_1\cdots j_q}|$). $R^n_+(R^n_{++})$
denotes the cone $\{x \in R^n\ | x_i\geq(>) 0, i=1,\ldots,n\}$. The
symbol $A \geq(>,\leq, <) B$ means that $a_{ij} \geq(>,\leq, <)
b_{ij}$ for every $i,j$ and it is the same for rectangular tensors.

\section{Preliminaries}
\setcounter{equation}{0} \setcounter{Assumption}{0}
\setcounter{Theorem}{0} \setcounter{Proposition}{0}
\setcounter{Corollary}{0} \setcounter{Lemma}{0}
\setcounter{Definition}{0} \setcounter{Remark}{0}
\setcounter{Algorithm}{0}\hspace*{\parindent} First we recall the
definition of tensor: a tensor is a multidimensional array, and a
real order m dimension n tensor $\mathcal{A}$ consists of $n^m$ real
entries:
\begin{displaymath}
A_{i_1\cdots i_m} \in R,
\end{displaymath}
where $i_j=1,\cdots,n$ for $j=1,\cdots,m$. If a number
$\lambda$ and a nonzero vector x are solutions of the following
homogeneous polynomial equations:
\begin{displaymath}
 \mathcal{A}x^{m-1} = \lambda x^{[m-1]},
\end{displaymath}
then $\lambda$ is called the eigenvalue of $\mathcal{A}$ and $x$ the
eigenvector of $\mathcal{A}$ associated with $\lambda$, where
$\mathcal{A}x^{m-1}$ and $x^{[m-1]}$ are vectors, whose $i$th
component are
\begin{eqnarray*}
(\mathcal{A}x^{m-1})_i &=& \sum^n_{i_2,\cdots,i_m=1}a_{ii_2\cdots i_m}x_{i_2}\cdots x_{i_m}\\
(x^{[m-1]})_i &=& x^{m-1}_i,
\end{eqnarray*}
respectively.  This definition was introduce by Qi \cite{s4} where
he supposed that $\mathcal{A}$ is an order $m$ dimension $n$
symmetric tensor and $m$ is even. Independently, Lim \cite{s2} gave
such a definition but restricted $x$ to be a real vector and
$\lambda$ to be a real number. Here we use the definition given in
\cite{s1}.

The Perron-Frobenius theorem for nonnegative tensors is related to
measuring high-order connectivity in linked objects and hypergraphs,
see \cite{s5,s6}. Let we recall the Perron-Frobenius theorem for
nonnegative tensors given in \cite{s1}:
\begin{Theorem}(see theorem 1.3 of \cite{s1}) \label{th:99}
If $\mathcal{A}$ is a nonnegative tensor of order m dimension n,
then there exists $\lambda_0 \geq 0$ and a nonnegative vector $x_0
\neq 0$ such that
\begin{equation} \label{eq:55555}
\mathcal{A}x^{m-1} = \lambda_0 x^{[m-1]}_0.
\end{equation}
\end{Theorem}
\begin{Theorem}(see theorem 1.4 of \cite{s1}) \label{th:100}
If $\mathcal{A}$ is an irreducible nonnegative tensor of order m
dimension n, then the pair $(\lambda_0, x_0)$ in equation
(\ref{eq:55555}) satisfy:
\begin{enumerate}
\renewcommand{\labelenumi}{$($\arabic{enumi}$)$}
\item $\lambda_0 > 0$ is an eigenvalue.
\item $x_0 > 0$, i.e. all components of $x_0$ are positive.
\item If $\lambda$ is an eigenvalue with nonnegative eigenvector,
then $\lambda = \lambda_0$. Moreover, the nonnegative eigenvector is
unique up to a multiplicative constant.
\item If $\lambda$ is an eigenvalue of $\mathcal{A}$, then $|\lambda| \leq
\lambda_0$.
\end{enumerate}
\end{Theorem}
And the reducibility of tensor is defined as follow:
\begin{Definition}(Reducibility, see definition 2.1 of \cite{s1}) A tensor $\mathcal {C} = (c_{i_1}\cdots
c_{i_m})$ of order m dimension n is called reducible, if there
exists a nonempty proper index subset $I \subset \{1,\cdots,n\}$
such that
\begin{displaymath}
c_{i_1\cdots i_m} = 0, \qquad \forall i_1 \, \in \, I, \quad \forall
i_2,\cdots,i_m \,\not\in \,I.
\end{displaymath}
\end{Definition}
If $\mathcal{C}$ is not reducible, then we call $\mathcal{C}$
irreducible. In \cite{s11}, Yang and Yang prove that for any
nonnegative tensor, the spectral radius is the largest eigenvalue of
it, which is an enhancement of theorem \ref{th:99}:
\begin{Theorem}(See theorem 2.3 of \cite{s11}) \label{th:101}
If $\mathcal{A}$ is a nonnegative tensor of order m dimension n,
then $\rho(\mathcal{A})$ is an eigenvalue with a nonnegative
eigenvector $y \in R^n_+$ corresponding to it.
\end{Theorem}
\begin{Definition}
The spectral radius of tensor $\mathcal{A}$ is defined as
\begin{displaymath}
\rho(\mathcal{A}) = \max\{|\lambda|: \lambda \textrm{ is an eigenvalue of }\mathcal{A}\}.
\end{displaymath}
\end{Definition}
For positive tensors, the following theorem holds:
\begin{Theorem}(See theorem 2.4 of \cite{s11}) \label{th:102}
Let $\mathcal{A}$ be a positive order m dimension n tensor, if
$\lambda$  is an eigenvalue of $\mathcal{A}$ except
$\rho(\mathcal{A})$, then $\rho(\mathcal{A}) > |\lambda|$.
\end{Theorem}

\section{Algorithm}
\setcounter{equation}{0} \setcounter{Assumption}{0}
\setcounter{Theorem}{0} \setcounter{Proposition}{0}
\setcounter{Corollary}{0} \setcounter{Lemma}{0}
\setcounter{Definition}{0} \setcounter{Remark}{0}
\setcounter{Condition}{0} \setcounter{Remark}{0}
\setcounter{Algorithm}{0}\hspace*{\parindent} Before presenting our
algorithm, we give the definition of diagonal similar tensors, which
was first used by Yang \emph{et al} \cite{s11}:
\begin{Definition}(diagonal similar tensors)
Let $\mathcal{A}=(a_{i_1\cdots i_m}), \mathcal{B}=(b_{i_1\cdots
i_m})$ be two order $m$ dimension $n$ tensors, if there is a
nonsingular diagonal matrix $D=(d_{ij})$, such that
\begin{displaymath}
\mathcal{A} = \mathcal{B}\cdot D^{-(m-1)} \cdot
\overbrace{D\cdot\cdots\cdot D}^{m-1},
\end{displaymath}
where
\begin{displaymath}
a_{i_1i_2\cdots i_m} = d^{-(m-1)}_{i_1,i_1}b_{i_1i_2\cdots
i_m}d_{i_2,i_2}\cdots d_{i_m,i_m},\,\,\, i_1,\cdots,i_m \, \in \,
\{1,\cdots,n\}.
\end{displaymath}
\end{Definition}
On the diagonal similar tensors we have the following proposition:
\begin{Proposition} \label{pro:2}
Let $\mathcal{A},\mathcal{B}, D$ be defined as above, if
$\mathcal{A},\mathcal{B}$ have eigenvalues, then they have the same
eigenvalues, i.e, if $\lambda$ be an eigenvalue of $\mathcal{A}$
with corresponding eigenvector $x \neq 0$, then $\lambda$ is also an
eigenvalue of $\mathcal{B}$ with corresponding eigenvector $Dx$; if
$\rho$ is an eigenvalue of $\mathcal{B}$ with corresponding
eigenvector $y \neq 0$, then $\rho$ is also an eigenvalue of
$\mathcal{A}$ with corresponding eigenvector $D^{-1}y$.
\end{Proposition}
\textbf{Proof.}
\begin{eqnarray*}
(\mathcal{B}\cdot(Dx)^{m-1})_i &=&
\sum_{i_2,\cdots,i_m=1}b_{i,i_2\cdots i_m}d_{i_2,i_2}\cdots
d_{i_m,i_m}x_{i_2}\cdots x_{i_m} \\
&=& d_{i,i}^{m-1}\sum_{i_2,\cdots,i_m=1}a_{i,i_2\cdots
i_m}x_{i_2}\cdots x_{i_m}\\
&=& \lambda d_{i,i}^{m-1}x^{m-1}_i = \lambda (Dx)_i,\,\,\,
i=1,\cdots,n.
\end{eqnarray*}
The proof of $\mathcal{A}(D^{-1}y)^{m-1}=\rho (D^{-1}y)^{[m-1]}$ is
the same. \ep

One easily gets following estimation:
\begin{Lemma} \label{lem:0.1}(See lemma 5.6 of \cite{s11})
Let $\mathcal{A} \geq 0$ be an order $m$ dimension $n$ tensor. Denote $\rho(\mathcal{A})$ the spectral radius of $\mathcal{A}$.
Then
\begin{displaymath} \min_{1\leq i \leq n}\sum_{i_2,\cdots,
i_m=1}a_{i,i_2\cdots i_m} \leq \rho(\mathcal{A}) \leq \max_{1\leq i
\leq n}\sum_{i_2,\cdots, i_m=1}a_{i, i_2\cdots i_m}.
\end{displaymath}
\end{Lemma}
Let $R_i = \sum_{i_2,\cdots,i_m=1}a_{i,i_2\cdots i_m}$. If $R_i
\equiv C$, where $C$ is a constant, we have the
following proposition:
\begin{Lemma} \label{lem:0.3} (See lemma 5.5 of \cite{s11})
Suppose $\mathcal{A} \geq 0$ be an order $m$ dimension $n$ tensor,
if $R_i=\sum_{i_2, \dots, i_m=1}a_{i,i_2\cdots i_m} \equiv
C(constant)$ for $i=1,\cdots,n$, then
\begin{displaymath}
\rho(\mathcal{A}) = C.
\end{displaymath}
\end{Lemma}
Given an order m dimension n  nonnegative irreducible tensor
$\mathcal{B}$, we calculate its spectral radius as follow:
\begin{Algorithm} \label{alg:1}
\textbf{step 1}. Let $\mathcal{A} = \mathcal{B} + \mathcal{I}$.
Denote $\mathcal{A}^{(0)} = \mathcal{A}$, set $k\,:= 0$, compute
\begin{displaymath}
R^{(0)}_i = \sum_{i_2,\cdots,i_m=1}a_{i,i_2,\cdots,i_m},
\end{displaymath}

\begin{displaymath}
R^{(0)} = \max_{1\leq i\leq n}R^{(0)}_i,\, r^{(0)} = \min_{1\leq i
\leq n}R^{(0)}_i.
\end{displaymath}
if $R^{(0)} = r^{(0)} $, let $v(\mathcal{A}) = e = (1,\cdots,1)^T$, goto step 3, else goto step2. \\
\textbf{step 2}. Compute
\begin{displaymath}
A^{(k+1)} = A^{(k)} \cdot D(k)^{-(m-1)}\overbrace{D(k)\cdot \cdots
D(k)}^{m-1},
\end{displaymath}
where
\begin{displaymath}
D(k) = diag((R^k_1)^{\frac{1}{m-1}}, \cdots,
(R^k_n)^{\frac{1}{m-1}})
\end{displaymath}
and
\begin{displaymath}
a^{(k+1)}_{i,i_2,\cdots,i_m} =
(R^{(k)})^{-1}a^{(k)}_{i,i_2,\cdots,i_m}(R^{(k)}_{i_2})^{\frac{1}{m-1}}\cdots
(R^{(k)}_{i_m})^{\frac{1}{m-1}},
\end{displaymath}
compute
\begin{displaymath}
R^{(k+1)}_i = \sum_{i_2,\cdots,i_m=1}a^{(k+1)}_{i,i_2\cdots i_m}
\end{displaymath}
\begin{displaymath}
R^{(k+1)} = \max_{1\leq i\leq n}R^{(k+1)}_i,\, r^{(k+1)} =
\min_{1\leq i\leq n}R^{(k+1)}_i
\end{displaymath}
\begin{displaymath}
v^{(k+1)}(\mathcal{A}) :=
diag((\prod^{k}_{j=0}\frac{R^{(j)}_1}{R^{(j)}})^{\frac{1}{m-1}},\cdots,(\prod^{k}_{j=0}\frac{R^{(j)}_n}{R^{(j)}})^{\frac{1}{m-1}})\cdot
((\frac{R^{(k+1)}_1}{R^{(k+1)}})^{\frac{1}{m-1}},\cdots,(\frac{R^{(k+1)}_n}{R^{(k+1)}})^{\frac{1}{m-1}})^T
\end{displaymath}
if $R^{(k+1)} = r^{(k+1)}$, goto step 3, else loop step 2.\\
\textbf{step 3}. Output $\rho(\mathcal{B})= R^{(k)} - 1$,
$v(\mathcal{A})$ as the eigenvalue and eigenvector of $\mathcal{B}$.
\end{Algorithm}
\textbf{Remark.} When $m=2$, this algorithm reduces to the smoothing method in
\cite{s14}, and $R^{(k)}_i$ is the sum of the $i$th row of matrix
$A^{(k)}$.

It is easy to notice that $\mathcal{A}$ satisfies the following
condition:
\begin{Condition}  \label{cdt:1}
\begin{displaymath}
R_i = \sum_{i_2,\cdots,i_m}a_{i,i_2\cdots i_m} > 0, \,\,\, i =
1,\cdots,n.
\end{displaymath}
\end{Condition}
Hence algorithm \ref{alg:1} is well-defined. Moreover, we have the
following theorem to ensure that $R^{(k)} - r^{(k)}$ is
nonincreasing as $k \rightarrow \infty$ (and in fact we can prove it
decreases strictly):
\begin{Theorem} \label{th:0.1} Under Condition.1,
\begin{displaymath}
r \leq r^{(1)} \leq \cdots \leq r^{(k)} \leq \cdots \leq
\rho(\mathcal{A}) \leq \cdots \leq R^{(k)} \leq \cdots \leq R.
\end{displaymath}
\end{Theorem}

Under the assumption of irreducibility of $\mathcal{B}$, we have
\begin{Theorem} \label{th:0.2} If $\mathcal{A} = \mathcal{B} + \mathcal{I}$ , then
\begin{displaymath}
\lim_{k \rightarrow \infty}r^{(k)} = \lim_{k \rightarrow
\infty}R^{(k)} = \rho(\mathcal{A}).
\end{displaymath}
\end{Theorem}
This theorem shows that the algorithm can find the spectral radius.
We will prove these theorems in the next section.

\section{Convergence analysis}
\setcounter{equation}{0} \setcounter{Assumption}{0}
\setcounter{Theorem}{0} \setcounter{Proposition}{0}
\setcounter{Corollary}{0} \setcounter{Lemma}{0}
\setcounter{Definition}{0} \setcounter{Remark}{0}
\setcounter{Algorithm}{0} \hspace*{\parindent} In this section, we
will prove theorems \ref{th:0.1} and \ref{th:0.2}. First we give a
lemma:
\begin{Lemma} \label{lem:0.2}
\begin{displaymath}
R^{(k-1)} \geq R^{(k)}, r^{(k-1)} \leq r^{(k)}, k=1,2,\cdots,
\end{displaymath}
where $R^{(k)}, r^{(k)}$ are defined in algorithm \ref{alg:1}.
\end{Lemma}
\textbf{Proof.} Without loss of generality we suppose that $R^{(k)}
= R^{(k)}_s$, $r^{(k)} = R^{(k)}_t$, $s,t \in \{1,2,\cdots,n\}$. We
have
\begin{eqnarray*}
R^{(k)} &=& \sum_{i_2,\cdots,i_m=1}a^{(k-1)}_{s,i_2\cdots
i_m}\frac{(R^{(k-1)}_{i_2})^{\frac{1}{m-1}}\cdots
(R^{(k-1)}_{i_m})^{\frac{1}{m-1}}}{R^{(k-1)}_s}
\\
&\leq&
\sum_{i_2,\cdots,i_m=1}a_{s,i_2\cdots i_m}\frac{R^{(k-1)}}{R^{(k-1)}_s}
\\&=& R^{(k-1)} \\
r^{(k)} &=& \sum_{i_2,\cdots,i_m=1}a^{(k-1)}_{t,i_2\cdots
i_m}\frac{(R^{(k-1)}_{i_2})^{\frac{1}{m-1}}\cdots
(R^{(k-1)}_{i_m})^{\frac{1}{m-1}}}{R^{k-1}_t} \\
&\geq&
\sum_{i_2,\cdots,i_m=1}a_{s,i_2\cdots i_m}\frac{r^{(k-1)}}{R^{k-1}_t}
\\&=& r^{(k-1)}.
\end{eqnarray*}\ep

\textbf{Proof of theorem \ref{th:0.1}:} From lemma \ref{lem:0.1} and
\ref{lem:0.2} we can easily see that
theorem \ref{th:0.1} holds. \ep\\

The following proposition shows that the sequence $\{R^{(k)} -
r^{(k)}\}$ is decreasing strictly:

\begin{Proposition} \label{prop:2222}
Let $\mathcal{A}$ be defined in algorithm \ref{alg:1}, $i \in I =
\{1,\cdots,n\}$, and without loss of generality suppose that
$R^{(k)} = R^{(k)}_s, r^{(k)} = R^{(k)}_t$; let $J = J(k) =
\{(i_2,\cdots,i_m)| \frac{a^{(k)}_{s,i_2\cdots i_m}}{R^{(k)}_s} \geq
\frac{a^{(k)}_{t,i_2\cdots i_m}}{R^{(k)}_t}\}$, $N
=\{(i_2,\cdots,i_m)| i_2,\cdots,i_m \in I\}$. Then we have
\begin{equation} \label{eq:4444}
R^{(k)} - r^{(k)} \leq (R^{(k-1)} - r^{(k-1)})[1-
\frac{1}{R^{(k-1)}}(\sum_{(i_2,\cdots,i_m) \in
N-J}a^{(k-1)}_{s,i_2\cdots i_m} + \sum_{(i_2,\cdots,i_m) \in
J}a^{(k-1)}_{t,i_2\cdots i_m})]
\end{equation}
\end{Proposition}
\textbf{Proof.} We only prove the case when $k=1$, i.e.
\begin{equation} \label{eq:1111}
R^{(1)} - r^{(1)} \leq (R - r)[1- \frac{1}{R}(\sum_{(i_2,\cdots,i_m)
\in N-J}a_{s,i_2\cdots i_m} + \sum_{(i_2,\cdots,i_m) \in
J}a_{t,i_2\cdots i_m})].
\end{equation}
For $k=2,3,\cdots$ the proof is the same.

If $R^{(1)} = r^{(1)}$, then (\ref{eq:1111}) holds easily, so we
suppose that $R^{(1)} > r^{(1)}$.
\begin{eqnarray} \label{eq:2222}
R^{(1)} - r^{(1)} &=& R^{(1)}_s - R^{(1)}_t \nonumber\\
&=& \sum_{i_2,\cdots,i_m=1}^n\frac{ a_{s,i_2\cdots
i_m}R_{i_2}^{\frac{1}{m-1}}\cdots R_{i_m}^{\frac{1}{m-1}}}{R_s} -
\sum_{i_2,\cdots,i_m=1}^n\frac{ a_{t,i_2\cdots
i_m}R_{i_2}^{\frac{1}{m-1}}\cdots R_{i_m}^{\frac{1}{m-1}}}{R_t} \nonumber\\
&=& \sum_{i_2,\cdots,i_m=1}^n(\frac{ a_{s,i_2\cdots i_m}}{R_s} -
\frac{ a_{t,i_2\cdots i_m}}{R_t} )R_{i_2}^{\frac{1}{m-1}}\cdots
R_{i_m}^{\frac{1}{m-1}}
\end{eqnarray}
By definition of $J$ and $R^{(1)} >r^{(1)}$ we see that $J, N-J \neq
\emptyset$, respectively, thus we have
\begin{equation} \label{eq:3333}
\sum^n_{(i_2,\cdots,i_m) \in J}(\frac{a_{s,i_2\cdots i_m}}{R_s} -
\frac{a_{t,i_2\cdots i_m}}{R_t}) = -\sum^n_{(i_2,\cdots,i_m) \in
N-J}(\frac{a_{s,i_2\cdots i_m}}{R_s} - \frac{a_{t,i_2\cdots
i_m}}{R_t}).
\end{equation}
Combining (\ref{eq:2222}) and (\ref{eq:3333}) we have
\begin{eqnarray*}
R^{(1)} - r^{(1)} &=& \sum_{(i_2,\cdots,i_m) \in
J}(\frac{a_{s,i_2\cdots i_m}}{R_s} - \frac{a_{t,i_2\cdots
i_m}}{R_t})R^{\frac{1}{m-1}}_{i_2}\dots R^{\frac{1}{m-1}}_{i_m} \\
&&+\sum_{(i_2,\cdots,i_m) \in N-J}(\frac{a_{s,i_2\cdots i_m}}{R_s} -
\frac{a_{t,i_2\cdots i_m}}{R_t})R^{\frac{1}{m-1}}_{i_2}\dots
R^{\frac{1}{m-1}}_{i_m} \\
&\leq& R\sum_{(i_2,\cdots,i_m) \in J}(\frac{a_{s,i_2\cdots
i_m}}{R_s} - \frac{a_{t,i_2\cdots i_m}}{R_t}) +
r\sum_{(i_2,\cdots,i_m) \in N-J}(\frac{a_{s,i_2\cdots i_m}}{R_s} -
\frac{a_{t,i_2\cdots i_m}}{R_t}) \\
&=& R\sum_{(i_2,\cdots,i_m) \in J}(\frac{a_{s,i_2\cdots i_m}}{R_s} -
\frac{a_{t,i_2\cdots i_m}}{R_t}) - r\sum_{(i_2,\cdots,i_m) \in
J}(\frac{a_{s,i_2\cdots i_m}}{R_s} - \frac{a_{t,i_2\cdots
i_m}}{R_t}) \\
&=& (R-r)(\sum_{(i_2,\cdots,i_m) \in J}\frac{a_{s,i_2\cdots
i_m}}{R_s} - \sum_{(i_2,\cdots,i_m) \in J}\frac{a_{t,i_2\cdots
i_m}}{R_t}) \\
&=& (R-r)[1-(\sum_{(i_2,\cdots,i_m) \in N-J}\frac{a_{s,i_2\cdots
i_m}}{R_s} + \sum_{(i_2,\cdots,i_m) \in J}\frac{a_{t,i_2\cdots
i_m}}{R_t})]\\
&\leq& (R-r)[1-\frac{1}{R}(\sum_{(i_2,\cdots,i_m) \in
N-J}a_{s,i_2\cdots i_m} + \sum_{(i_2,\cdots,i_m) \in
J}a_{t,i_2\cdots i_m}).
\end{eqnarray*}
The proof is completed. \ep\\

From Proposition \ref{prop:2222} we see that $\{R^{(k)} - r^{(k)}\}$
is a nonnegative and strictly monotone decreasing sequence, so it
has limit, but it is not sufficient to show theorem \ref{th:0.2}.
Before prove it, we denote
$M^{(k)}_i=(\prod^{k}_{l=0}\frac{R^{(l)}_i}{R^{(l)}})^{\frac{1}{m-1}},i=1,2,\cdots,n$,
hence $v^{(k)}(\mathcal{A})$ defined in algorithm \ref{alg:1} is
$(M^{(k)}_1,\cdots, M^{(k)}_n)^T$ and $V^{(k)}(\mathcal{A})$ the
matrix representation of $v^{(k)}(\mathcal{A})$, i.e:
\begin{displaymath}
V^{(k)}(\mathcal{A}) = diag(M^{(k)}_1,\cdots,M^{(k)}_n).
\end{displaymath}
Note that the limit of
$M^{(k)}_i=(\prod^{k}_{l=0}\frac{R^{(l)}_i}{R^{(l)}})^{\frac{1}{m-1}}$
exists. Denote it by $M^*_i, i=1,\cdots,n$, and $v^*(\mathcal{A})$,
$V^*(\mathcal{A})$ the limit of $v^{(k)}(\mathcal{A})$,
$V^{(k)}(\mathcal{A})$, respectively. We introduce a notation:
\begin{eqnarray*}
&&V(\mathcal{A})^{m-1}\cdot \mathcal{A}\cdot V(\mathcal{A})^{-(m-1)}
\cdot \overbrace{V(\mathcal{A}) \cdot\cdots\cdot
V(\mathcal{A})}^{m-1}\\
&=& \mathcal{A}\cdot (V(\mathcal{A})^{m-1}V(\mathcal{A})^{-(m-1)})
\cdot \overbrace{V(\mathcal{A}) \cdot\cdots\cdot
V(\mathcal{A})}^{m-1}\\
&=& \mathcal{A}\cdot I \cdot \overbrace{V(\mathcal{A})
\cdot\cdots\cdot V(\mathcal{A})}^{m-1}.
\end{eqnarray*}
Where the superscript '$(k)$' is ignored. Then we have
\begin{eqnarray*}
\mathcal{A}v^{(k)}(\mathcal{A})^{m-1} &=&
\mathcal{A}(V^{(k)}(\mathcal{A})e)^{m-1} \\
&=& \mathcal{A}\cdot I\cdot
\overbrace{V^{(k)}(\mathcal{A})\cdot\cdots\cdot
V^{(k)}(\mathcal{A})}^{m-1}\cdot
e^{m-1}\\
&=& V^{(k)}(\mathcal{A})^{m-1}\cdot (\mathcal{A}\cdot
V^{(k)}(\mathcal{A})^{-(m-1)} \cdot
\overbrace{V^{(k)}(\mathcal{A})\cdot\cdots \cdot
V^{(k)}(\mathcal{A})}^{m-1})e^{m-1}\\
&=& V^{(k)}(\mathcal{A})^{m-1}\cdot \mathcal{A}^{(k+1)}e^{m-1}\\
&=& V^{(k)}(\mathcal{A})^{m-1}\cdot (R^{(k+1)}_1,\cdots,R^{k+1}_n)^T\\
&\geq& r^{(k+1)}v^{(k)}(\mathcal{A})^{[m-1]}.
\end{eqnarray*}
Letting $k \rightarrow \infty$ we have
\begin{equation} \label{eq:6666}
\mathcal{A}v^{*}(\mathcal{A})^{m-1} \geq
r^{*}v^{*}(\mathcal{A})^{[m-1]}.
\end{equation}
By theorem \ref{th:0.1},  $r^*$ exists, where $r^*$ is the limit of
$r^{(k)}$, and it satisfies the proposition:
\begin{Proposition} \label{pro:3}
$r^* \geq r^{(k)},k=1,2,\cdots.$
\end{Proposition}
 If (\ref{eq:6666})
is an equation, by theorem \ref{th:100}, $r^*$ is the spectral
radius. Suppose not, denote
\begin{eqnarray*}
x_0(v^{*}(\mathcal{A})) = (\mathcal{A}v^{*}(\mathcal{A})^{m-1})^{[\frac{1}{m-1}]},\\
y_0(v^{*}(\mathcal{A})) = (r^{*}v^{*}(\mathcal{A})^{[m-1]})^{[\frac{1}{m-1}]};\\
x_1(v^{*}(\mathcal{A})) = (\mathcal{A}x_0^{m-1})^{[\frac{1}{m-1}]},\\
y_1(v^{*}(\mathcal{A})) = (\mathcal{A}y_0^{m-1})^{[\frac{1}{m-1}]};\\
\cdots \\
x_{(n-1)}(v^{*}(\mathcal{A})) = (\mathcal{A}x_{(n-2)}^{m-1})^{[\frac{1}{m-1}]},\\
y_{(n-1)}(v^{*}(\mathcal{A})) =
(\mathcal{A}y_{(n-2)}^{m-1})^{[\frac{1}{m-1}]}.
\end{eqnarray*}
We will prove
\begin{equation} \label{eq:7777}
x_{(n-1)}(v^{*}(\mathcal{A})) > y_{(n-1)}(v^{*}(\mathcal{A}))
\end{equation}. Recall theorem 6.6 given in
\cite{s11}:
\begin{Theorem}(see theorem 6.6 of \cite{s11}) \label{th:4444}
Let $\mathcal{B} \geq 0$ be an order m dimension n tensor. Then
$\mathcal{B}$ is irreducible if and only if for all $x \in R^n_+$,
$x \neq 0$, let $x_0 = x$ and $x_{k+1} =
(\mathcal{B}+\mathcal{I})x_{k}^{m-1}$. Then $x_{n-1}
> 0$.
\end{Theorem}
We follow this theorem to prove a lemma:
\begin{Lemma} \label{lem:5555}
Let $\mathcal{B}$ be defined as above. Denote $\mathcal{A} =
\mathcal{B} + \mathcal{I}$. Suppose $x,y \in R^n_+$ and $x \geq y$.
let $x_0 = x,y_0=y$ and $x_{k+1} = \mathcal{A}x_{k}^{m-1}$, $y_{k+1}
= \mathcal{A}y_{k}^{m-1}$. Then $x_{(n-1)} > y_{(n-1)}$.
\end{Lemma}
\textbf{Proof.} It is easy to see that $x_{(n-1)} \geq y_{(n-1)}$.
Let $I_k = \{i | x_{k_i} = y_{k_i}\}$ and $M_k = |I_k|$, so does
$M_{k+1}$. All we need to do is to prove that $M_{k+1} < M_k$. First
we consider $x_{k+1_i} = (\mathcal{A}x_k^{m-1} + x_k^{[m-1]})_i$
where $i \not\in I_k$. We have $x_{k+1_i} > y_{k+1_i}$ because
$(\mathcal{A}x_k^{m-1})_i \geq (\mathcal{A}y_k^{m-1})_i $ and
$x_{k_i}
> y_{k_i}$. For all $i \in I_k$, $x_{k+1_i} = \sum_{i2,\cdots,i_m =
1}a_{ii2 \cdots i_m}x_{k_{i_2}}\cdots x_{k_{i_m}}$. We claim that at
least a $x_{k+1_i} > y_{k+1_i}$, $i \in I_k$. Suppose not, then
$\sum_{i2,\cdots,i_m = 1}a_{ii2 \cdots i_m}(x_{k_{i_2}}\cdots
x_{k_{i_m}} - y_{k_{i_2}}\cdots y_{k_{i_m}})= 0$, $i \in I_k$. It
means that $a_{ii_2\cdots i_m }=0$, $i_2,\cdots,i_m \not\in I_k$ and
$i \in I_k$, which contradicts to the irreducibility of
$\mathcal{A}$. Thus at least a $x_{k+1_i}>y_{k+1_i} $, $i \in I_k$,
which means that $M_{k+1} < M_k$. Repeat at most $n-1$ times, we
have $x_{n-1}
> y_{n-1}$. \ep

It is easy to see that this lemma holds when we replace $x_{k+1} =
\mathcal{A}x_{k}^{m-1}$, $y_{k+1} = \mathcal{A}y_{k}^{m-1}$ by
$x_{k+1} = (\mathcal{A}x_{k}^{m-1})^{[\frac{1}{m-1}]}$, $y_{k+1} =
(\mathcal{A}y_{k}^{m-1})^{[\frac{1}{m-1}]}$. By this lemma, if
(\ref{eq:6666}) is not an equation, then (\ref{eq:7777}) holds. By
the continuity of $x_{(n-1)}(\cdot) - y_{(n-1)}(\cdot)$, when $k$
sufficiently large, one has
\begin{equation}\label{eq:8888}
x_{(n-1)}(v^{k}(\mathcal{A}))
> y_{(n-1)}(v^{k}(\mathcal{A})),
\end{equation} and
\begin{eqnarray*}
x_0(v^{k}(\mathcal{A})) =
(\mathcal{A}v^{(k)}(\mathcal{A})^{m-1})^{[\frac{1}{m-1}]} =
(R^{(k+1)})^{\frac{1}{m-1}}v^{(k+1)}(\mathcal{A}),\\
x_1(v^{k}(\mathcal{A})) =
(R^{(k+1)}R^{(k+2)})^{\frac{1}{m-1}}v^{(k+2)}(\mathcal{A}),\\
\cdots \\
x_{(n-1)}(v^{k}(\mathcal{A})) =
(\prod^n_{l=1}R^{(k+l)})^{\frac{1}{m-1}}v^{(k+n)}(\mathcal{A}),\,\,and\\
y_{(n-1)}(v^{k}(\mathcal{A})) =
r^*(\prod^{n-1}_{l=1}R^{(k+l)})^{\frac{1}{m-1}}v^{(k+n-1)}(\mathcal{A}).
\end{eqnarray*}
Then (\ref{eq:8888}) means that
\begin{displaymath}
R^{(k+n)}_i > r^*,\,\,\,\,i=1,\cdots,n,
\end{displaymath}
especially $r^{(k+n)} > r*$ when $k$ sufficiently large, which
contradicts with proposition \ref{pro:3}. Thus
\begin{displaymath}
\mathcal{A}v^{*}(\mathcal{A})^{m-1}
=r^{*}v^{*}(\mathcal{A})^{[m-1]}.
\end{displaymath}
The same proof can apply to $\mathcal{A}v^{*}(\mathcal{A})^{m-1}
=R^{*}v^{*}(\mathcal{A})^{[m-1]}.$ Hence theorem \ref{th:0.2} holds
and $v^*(\mathcal{A})$ is the positive corresponding eigenvector of
$\rho(\mathcal{A})$.

\section{Numerical results}\hspace*{\parindent}
In this section, we first give numerical result on a 3-order
3-dimension nonnegative irreducible tensor; then we generate some
random tensors to test our algorithm. We use the termination 
condition given in \cite{s12}:\\
(1) $k \geq 100,$\\
(2) $R^{(k)} - r^{(k)} \leq 10^{(-7)}$.\\
\textbf{Example 1}. Consider the 3-order 3-dimensional tensor
\begin{displaymath}
\mathcal{B} = [B(1,:,:),B(2,:,:),B(3,:,:)],
\end{displaymath}
where
\begin{displaymath}
B(1,:,:) = \left(\begin{array}{ccc}
0 & 0 & 0 \\
0 & 3.72 & 0 \\
0 & 0 & 0
\end{array} \right),
\end{displaymath}
\begin{displaymath}
B(2,:,:) = \left(\begin{array}{ccc}
9.02 & 0 & 0 \\
0 & 0 & 0 \\
0 & 0 & 0
\end{array} \right),
\end{displaymath}
\begin{displaymath}
B(3,:,:) = \left(\begin{array}{ccc}
9.55 & 0 & 0 \\
0 & 0 & 0 \\
0 & 0 & 0
\end{array} \right).
\end{displaymath}
\textbf{Example 2.} We use some randomly generated tensors to test
algorithm \ref{alg:1}. Each entry of these tensors is between $0$
and $10$.

Table 1 is the numerical results of algorithm \ref{alg:1} for
example 1 where the tensor is $\mathcal{A} = \mathcal{B} +
\mathcal{I}$. From the result, we get $\rho(\mathcal{B}) = 6.79262 -
1 = 5.79262$ and the positive corresponding eigenvector is
$(0.46224, 0.57681, 0.593515)^T$. But if we directly apply the
algorithm to $\mathcal{B}$, it does not converges. Table 2 shows
some numerical results on randomly generated tensors of different
order and dimension. In our experiment, the
algorithm converges to the spectral radius for all the tensors although there may be some reducible tensors in these randomly generated tensors.\\\\

\begin{table}[!h]
\tabcolsep 5mm \caption{numerical results of algorithm \ref{alg:1}
for example 1}
\begin{center}
\begin{tabular}{rlrlrlrlrl}
\hline
\multicolumn{2}{c}{$k$}&\multicolumn{2}{c}{$r^{(k)}$}&\multicolumn{2}{c}{$R^{(k)}$}&\multicolumn{2}{c}{$R^{(k)}-r^{(k)}$}&\multicolumn{2}{c}{$0.5*(R^{(k)}+r^{(k)})$}
\\ \hline
                 &1&   &                  4.72   &                        10.55   & &
          5.83   &  &    &                      7.635\\
&2&   &                  5.24894   &                        8.89712
& &
               3.64818   &  &    &                      7.07303\\
&3&   &                  5.65898   &                        8.2097
& &
              2.55071   &  &    &                      6.93434\\
&4&   &                  5.96904   &                        7.7527
& &
              1.78366   &  &    &                      6.86087\\
&5&   &                  6.19911   &                        7.45402
& &
               1.25491   &  &    &                      6.82656\\
&6&   &                  6.36745   &                        7.25147
& &
               0.88402   &  &    &                      6.80946\\
               ...\\
                 &48&   &                  6.79262   &                        6.79262   & &
                3.83995e-007   &  &    &                      6.79262\\
&49&   &                  6.79262   &                        6.79262
& &
                2.70932e-007   &  &    &                      6.79262\\
&50&   &                  6.79262   &                        6.79262
& &
                1.9116e-007   &  &    &                      6.79262\\
&51&   &                  6.79262   &                        6.79262
& &
                1.34875e-007   &  &    &                      6.79262\\
&52&   &                  6.79262   &                        6.79262
& &
                9.51629e-008   &  &    &                      6.79262\\
                    \hline
\end{tabular}
\end{center}
\end{table}

\begin{table}[!h]
\tabcolsep 5mm \caption{numerical results of algorithm \ref{alg:1}
for some randomly generated tensors}
\begin{center}
\begin{tabular}{lllll}
\hline
\multicolumn{1}{c}{$(n,m)$}&\multicolumn{1}{c}{$k$}&\multicolumn{1}{c}{$\rho(\mathcal{A})$}&\multicolumn{1}{c}{$R^{(k)}-r^{(k)}$}&\multicolumn{1}{c}{$\|
\mathcal{A}v(\mathcal{A})^{m-1}-\rho(\mathcal{A})v(\mathcal{A})^{[m-1]}\|_{\infty}$}
\\ \hline
 (10,3)  & 7 &   444.247       & 9.56E-09  &  5.21E-09\\
 (5,3)   & 8 &   111.111       & 2.56E-08  &  1.52E-08\\
(20,3)   & 6 &   1791.51       & 2.15E-08  &  9.03E-09\\
(30,3)   & 6 &   4057.86       & 1.34E-09  &  6.32E-10\\
 (50,3)  & 5 &   11225.7       & 3.60E-08  &  1.81E-08\\
(100,3)  & 5 &   45013         & 2.57E-09  &  1.21E-09\\
 (20,4)  & 6 &   35983.8       & 2.28E-11  &  6.85E-11\\
  (5,4)  & 7 &   556.015       & 2.44E-10  &  1.58E-10\\
(10,4)   & 6 &   4494.69       & 1.92E-10  &  1.04E-10\\
 (15,4)  & 6 &   15144         & 1.78E-11  &  2.38E-11\\
  (30,4) & 5 &   121554        & 1.02E-10  &  1.75E-10\\
   (5,5) & 5 &   2765.93       & 1.83E-08  &  9.43E-09\\
 (10,5)  & 5 &   44913.9       & 7.23E-10  &  7.02E-10\\
 (15,5)  & 4 &   227923        & 7.56E-09  &  6.42E-09\\
 (20,5)  & 4 &   720407        & 2.95E-08  &  4.01E-08\\
(5,6)    & 4 &   14091.6       & 2.69E-08  &  1.11E-08\\
 (10,6)  & 4 &   450133        & 1.27E-08  &  2.16E-08\\
  (15,6) & 4 &   3.42E+06      & 5.09E-08  &  4.03E-07\\
                    \hline
\end{tabular}
\end{center}
\end{table}

%\mbox{}\\\mbox{}\\\mbox{}\\\mbox{}\\\mbox{}\\\mbox{}\\\mbox{}\\\mbox{}\\\mbox{}\\\mbox{}\\\mbox{}\\\mbox{}\\

\section{Conclusion and remarks}\hspace*{\parindent}
We give an algorithm to find the spectral radius of nonnegative
tensors. When the tensor is irreducible, our algorithm can assure to
find out the spectral radius and its corresponding positive
eigenvector. This result is better than the algorithm proposed by
Michael \emph{et al} \cite{s12}. We can also apply algorithm
\ref{alg:1} to $\mathcal{A}(\alpha) = \mathcal{B} +
\alpha\mathcal{I}$ where $\alpha$ is a positive number. The choice
of $\alpha$ will affect the convergence rate of the algorithm.  This needs further
research.

\end{document}